\documentclass[11pt]{article}
\usepackage{graphicx}
\usepackage{amsmath,amssymb}
\def\1{{\bf 1}}

\def\be{\begin{equation}}
\def\ee{\end{equation}}

\begin{document}
\noindent {\Large \bf The Hybrid Idea of (Energy Minimization) Optimization Methods Applied to Study Prion Protein Structures Focusing on the $\beta$2-$\alpha$2 Loop}\\

\bigskip

\noindent {\Large Jiapu Zhang$^\text{ab*}$}\\
{\small {\it
\noindent $^\text{a}$Molecular Model Discovery Laboratory, Department of Chemistry \& Biotechnology,  Faculty of Science, Engineering \& Technology, Swinburne University of Technology, Hawthorn Campus, Hawthorn, Victoria 3122, Australia; 

\noindent $^\text{b}$Graduate School of Sciences, Information Technology and Engineering \& Centre of Informatics and Applied Optimisation, Faculty of Science, The Federation University Australia, Mount Helen Campus, Mount Helen, Ballarat, Victoria 3353, Australia;

\noindent $^\text{*}$Tel: +61-3-9214 5596, +61-3-5327 6335; 
jiapuzhang@swin.edu.au, j.zhang@federation.edu.au\\
}}

\noindent {\bf Abstract:} 
 {\it In molecular mechanics, current generation potential energy functions provide a reasonably good compromise between accuracy and effectiveness. This paper firstly reviewed several most commonly used classical potential energy functions and their optimization methods used for energy minimization. To minimize a potential energy function, about 95\% efforts are spent on the Lennard-Jones potential of van der Waals interactions; we also give a detailed review on some effective computational optimization methods in the Cambridge Cluster Database to solve the problem of Lennard-Jones clusters. From the reviews, we found the hybrid idea of optimization methods is effective, necessary and efficient for solving the potential energy minimization problem and the Lennard-Jones clusters problem. An application to prion protein structures is then done by the hybrid idea. We focus on the $\beta$2-$\alpha$2 loop of prion protein structures, and we found (i) the species that has the clearly and highly ordered $\beta$2-$\alpha$2 loop usually owns a 3$_{10}$-helix in this loop, (ii) a ``$\pi$-circle" Y128--F175--Y218--Y163--F175--Y169--R164--Y128(--Y162) is around the $\beta$2-$\alpha$2 loop.\\ 
}

\noindent {\bf Key words:} {\it Hybrid idea; computational optimization methods; energy minimization; potential energy; Lennard-Jones clusters; application to protein structures; prion proteins; the $\beta$2-$\alpha$2 loop.}\\

\section{Introduction}
In molecular mechanics, current potential energy functions provide a reasonably good accuracy to structural data obtained from X-ray crystallography and nuclear magnetic resonance (NMR), and dynamic data obtained from spectroscopy and inelastic neutron scattering and thermodynamic data. Currently, AMBER, CHARMM, GROMOS and OPLS/AMBER are among the most commonly used classical potential energy functions \cite{abraham_etal_GROMACS_2014, amber14, locatellis2008, paquetv2015, stote_etal1999}. The energy, $E$, is a function of the atomic positions, $R$, of all the atoms in the system these are usually expressed in term of Cartesian coordinates. The value of the energy is calculated as a sum of bonded (or internal) terms $E_{bonded}$, which describe the bonds, angles and bond rotations in a macromolecule, and a sum of non-bonded (or external) long-range terms $E_{non-bonded}$ \cite{abraham_etal_GROMACS_2014, amber14, stote_etal1999}:
{\small
\begin{eqnarray}\label{potential_energy_formula}
E_{potential}     & = & E_{bonded} + E_{non-bonded}  \nonumber\\
                  & = & (E_{bond-stretch} + E_{angle-bend} + E_{rotate-along-bond}  \nonumber\\
                  &   & (+ E_{Urey-Bradley} + E_{improper} + U_{CMAP} ) ) \nonumber\\
                  & & + (E_{van-der-Waals} + E_{electrostatic} + E_{hydrogen-bonds}). 
\end{eqnarray}            
 }
\noindent For example, for AMBER and CHARMM force fields \cite{amber14, stote_etal1999} respectively, the potential energy functions are \cite{amber14}:
\begin{eqnarray}\label{AMBER_potential_energy_formula}
E_{AMBER} &=&  \sum_{bonds} k_b(b-b_0)^2 + \sum_{angles} k_{\theta}(\theta - \theta_{eq})^2 +\sum_{dihedrals} (V_n/2) [1+\cos (n\phi -\gamma )] \nonumber \\ 
          & &+ \sum_{i=1}^{N-1} \sum_{j=i+1}^N \left [ \frac{A_{ij}}{R_{ij}^{12}} -\frac{B_{ij}}{R_{ij}^6} + \frac{q_iq_j}{\varepsilon R_{ij} } \right ],
\end{eqnarray}          
\begin{eqnarray}\label{CHARMM_potential_energy_formula}
E_{CHARMM} &=&  \sum_{bonds} k_b(b-b_0)^2 + \sum_{angles} k_{\theta} (\theta - \theta_0)^2 +\sum_{dihedrals} k_{\phi} [1+\cos (n\phi -\delta )] \nonumber \\
           & & +\sum_{Urey-Bradley} k_u (u-u_0)^2 +\sum_{impropers} k(\omega -\omega_0)^2 +\sum_{\phi, \psi} V_{CMAP}(\phi, \psi ) \nonumber \\   
           & &+ \sum_{nonbonded} \left ( \varepsilon \left [ \left ( \frac{R_{min_{ij}}}{R_{ij}} \right )^{12}
                                                   - \left ( \frac{R_{min_{ij}}}{R_{ij}} \right )^{6} 
                                            \right ] 
                                                   + \frac{q_iq_j}{\varepsilon R_{ij}} \right ),
\end{eqnarray} 
where $b, \theta, \phi, R_{ij}, u, \omega, \psi$ are basic variables ($b$ is the bond length of two atoms, $\theta$ is the angle of three atoms, $\phi$ is the dihedral angle of four atoms, $R_{ij}$ is the distance between atoms $i$ and $j$), and all other mathematical symbols are constant parameters specified in various force fields respectively. This paper will discuss how to effectively and efficiently use computational optimization methods to solve the minimization problem of the potential energy in Eq. (\ref{potential_energy_formula}), i.e. 
\begin{eqnarray}\label{minimize_potential_energy}
min \quad E_{potential}.
\end{eqnarray}

\vskip 1cm

Firstly, for Eq. (\ref{minimize_potential_energy}), we consider why we should perform energy minimization (EM). There are a number of reasons:
\begin{enumerate}
\item[] (i) To remove nonphysical (or bad) contacts / interactions. For example, when a structure that has been solved via X-ray crystallography, in the X-ray crystallization process, the protein has to be crystallized so that the position of its constituent atoms may be distorted from their natural position and contacts with neighbors in the crystal can cause changes from the in vitro structure; consequently, bond lengths and angles may be distorted and steric clashes between atoms may occur. Missing coordinates obtained from the internal coordinate facility may be far from optimal. Additionally, when two sets of coordinates are merged (e.g., when a protein is put inside a water box) it is possible that there are steric clashes / overlap presented in the resulting coordinate set (www.charmmtutorial.org).
\item[] (ii) In molecular dynamics (MD) simulations, if a starting configuration is very far from equilibrium, the forces may be excessively large and the MD simulation may fail \cite{abraham_etal_GROMACS_2014}. 
\item[] (iii) To remove all kinetic energy from the system and to reduce the thermal noise in the structures and potential energies \cite{abraham_etal_GROMACS_2014}.
\item[] (iv) Re-minimize is needed if the system is under different conditions. For example, in quantum mechanics / molecular mechanics (QM/MM) one part of the system is modeled in QM while the rest is modeled in MD, re-minimize the system with a new condition is needed (www.charmmtutorial.org).
\end{enumerate}
To perform EM is to make the system reaching to a equilibration state. EM of Eq. (\ref{minimize_potential_energy}) can be challenging, as there are many local minima that optimization algorithms might get stuck in without finding the global minima - in most cases, this is what will actually happen. Thus, how much and how far we should minimize should be well considered. Over-minimization can lead to unphysical ``freezing" of the structure and move too much from its original conformation; if not minimized enough and exactly, for example, the normal mode calculation cannot arrive at the bottom of its harmonic well. However, in MD, because the output of minimization is to be used for dynamics, it is not necessary for the optimization to be fully converged but a few hundreds or tens of local optimization search are good and kind enough. To make enough local optimization, usually, after we put the protein into a solvent (e.g. waters), first we restrain the protein by holding the solute fixed with strong force and only optimize the solvent, next holding the solute heavy atoms only, and then holding the CA atoms only, and lastly remove all restraints and optimize the whole system.   

\vskip 1cm

Secondly, for Eq. (\ref{minimize_potential_energy}), we consider what optimization algorithms we should use. In packages of \cite{abraham_etal_GROMACS_2014, bhandarkar_etal_NAMD_2012, amber14} etc, the following three local search optimization methods have been used.
\begin{enumerate}
\item[] (i) SD (steepest descent) method is based on the observation that if the real-valued function $E(x)$ is defined and differentiable in a neighborhood of a point $x_0$ then $E(x)$ decreases fastest if one goes from $x_0$ in the direction of the negative gradient of $E(x)$ at $x_0$. SD method is the simplest algorithm, it simply moves the coordinates in the negative direction of the gradient (hence in the direction of the force - the force is the (negative) derivative of the potential), without consideration of build ups in previous steps - this is the fastest direction making the potential energy decrease. SD is robust and easy to implement. But SD is not the most efficient especially when closer to minimum and in the vicinity of the local minimum. This is to say, SD does not generally converge to a local minimum, but it can rapidly improve the conformation when system is far from a minimum - quickly remove bad contacts and clashes. 
\item[] (ii) Conjugate gradient (CG) method is a method adds an orthogonal vector to the current direction of optimization search and then moves them in another direction nearly perpendicular to this vector. CG method is fast-converging and uses gradient information from previous steps. CG brings you very close to the local minimum, but performs worse far away from the minimum. CG is slower than SD in the early stages of minimization, but becomes more efficient closer to the energy minimum. In GROMACS CG cannot be used with constraints and in this case SD is efficient enough. When the forces are truncated according to the tangent direction, making it impossible to define a Lagrangian, CG method cannot be used to find the EM path.
\item[] (iii) L-BFGS method is a Quasi-Newton method that approximates the reverse of Hessian matrix $[\nabla^2 E(x)]^{-1}$ of $E(x)$ for the Newton method search direction $-[\nabla^2 E(x)]^{-1} \nabla E(x)$. L-BFGS method is mostly comparable to CG method, but in some cases converges 2$\sim$3 times faster with super-linear convergent rate (because it requires significantly fewer line search steps than Polak-Ribiere CG). L-BFGS of Nocedal approximates the inverse Hessian by a fixed number of corrections from previous steps. In practice L-BFGS converges faster than CG.
\item[] (iv) The combination of CG and LBFGS, so-called lbfgs-TNCG-BFGS method is a preconditioned truncated Newton CG method, it requires fewer minimization steps than Polak-Ribiere CG method and L-BFGS method, but L-BFGS can sometimes be faster in the terms of total CPU times. 
\end{enumerate}
If a global optimization is required, approaches such as simulated annealing (SA), parallel tempering method (super SA, also called replica exchange \cite{zhang2011c}), Metropolis algorithms and other Monte Carlo methods, Simplex method, Nudged Elastic Band method, different deterministic methods of discrete or continuous optimization etc may be utilized. The main idea of SA refinement is to heat up the system such that the molecule of interest has enough energy to explore a wide range of configurational space and get over local optimal energy barriers.  Relatively large structural rearrangements are permitted at these high temperatures. As the temperature is cooled gradually, the structural changes proceed in smaller steps, continuing to descend toward the global energy minimum.

\vskip 1cm

For solving Eq. (\ref{minimize_potential_energy}), without considering $E_{hydrogen-bonds}=\sum_{i=1}^{N-1} \sum_{j=i+1}^N \left [ \frac{C_{ij}}{R_{ij}^{12}} -\frac{D_{ij}}{R_{ij}^{10}} \right ]$, about 95\% of the CPU time of calculations is spent at 
\begin{eqnarray}\label{minimize_LJ_potential_energy}
min \quad E_{van-der-Waals}=\sum_{i=1}^{N-1} \sum_{j=i+1}^N \left [ \frac{A_{ij}}{R_{ij}^{12}} -\frac{B_{ij}}{R_{ij}^6} \right ],
\end{eqnarray}
where $C_{ij}, D_{ij}$ are constants. In \cite{locatellis2008}, this problem is also called Lennard-Jones (LJ) Atomic Cluster Optimization problem (where within the field of atomic clusters only nonbonded interactions are accounted for and particles are considered to be charge-free; e.g. real clusters of metals like gold, silver, and nickel). It is very necessary to up to date review some effective and efficient computational methods for solving Eq. (\ref{minimize_LJ_potential_energy}). There are numerous algorithms to solve Eq. (\ref{minimize_LJ_potential_energy}); here we just list the ones in  The Cambridge Energy Landscape Database (http://doye.chem.ox.ac.uk/jon/structures/LJ.html) which can obtain the best global structures: 
\begin{enumerate}
\item[$-$] Hoare  and  Pal's work \cite{hoarep1971a, hoarep1971b, hoarep1972} may  be the early  most successful results on LJ problem. The  idea is using build-up techniques to construct the initial solutions which are expected to represent low energy states, and using those initial solutions as starting points for a local search  method to relax to the optimal solution \cite{hoarep1971b}. The starting seed is the regular unit tetrahedron with atoms at the vertexes, the obvious global optimal solution for $N=4$. Beginning with this tetrahedron, Hoare and Pal (1971, 1972) added one atom at a time to construct a sequence of polytetrahedral  structures and at last got  good results up to $N=66$ \cite{hoarep1971a, hoarep1971b, hoarep1972}. For example, for $N=5$ its globally optimal trigonal bi-pyramid (bi-tetrahedron) structure is gotten by adding  an atom at the tetrahedral capping  position over a triangular face; following the bi-tetrahedron structure, the optimal structure of $N=6$ is tri-tetrahedron (another known optimal structure for $N=6$ is octahedron (using tetrahedral capping over triangular faces and half-octahedral capping over square faces), which is not a polytetrahedron); for $N=7$ its best structure constructed is the pentagonal bi-pyramid, a structure with a five-fold axis of symmetry.  Many computer science data structure procedures such as greedy forward growth operator and reverse greedy operator can make the build-up technique work well. The application of methods of studying noncrystalline clusters to the study of ``spherical" face centred cubic (fcc) microcrystallites was described in \cite{hoarep1972}. In \cite{hoarep1971a} the chief geometrical features of the clustering of small numbers of interacting particles were described.
\item[$-$] The  data structure of Northby \cite{northby1987} in finding the good starting solution  is the lattice based structure. The lattice structures consist of an icosahedral core and particular combinations of surface lattice points. A class of icosahedral packings was by constructed in \cite{mackay1962} adding successively larger icosahedral shells in layers around  a core central atom; this icosahedral  lattice can  be described  as 20 slightly flattened tetrahedrally shaped  fcc units with 12 vertices on a sphere centered at the core atom.  Atoms within each triangular face are placed in staggered rows in a two dimensional hexagonal close-packed arrangement.  Each atom in the interior of a face in a given shell is a tetrahedral capping  position relative to three atoms in the underlying  shell.  Northby (1987) relaxed  the structure of \cite{mackay1962} to get  his IC and  FC  multilayer  icosahedral  lattice  structures \cite{northby1987}. The IC lattice  can be referred  to the FORTRAN code in \cite{xue1994b}; it consists of all those  sites which will comprise the outer shell of the next complete Mackay \cite{mackay1962} icosahedron.  FC lattice is a slight modification of IC lattice in that its outer shell maintains icosahedral symmetry and consists of points at the icosahedral  vertices and the stacking fault positions of the outer IC shell. Basing on the IC and FC lattices, Northy (1987) gave his algorithm first finding a set of lattice local minimizers and then relaxing those lattice minimizers by performing  continuous minimization starting with those lattice minimizers \cite{northby1987}. The algorithm was summarized  as Algorithm 1 and Algorithm 2 of \cite{xue1994b}.
\item[$-$] The great majority of the best known solutions of Northy \cite{northby1987} are icosahedral in character. The hybridization of global search  and  local search  methods, usually,  is more effective to solve the large scale problem than the global search method or local search method working  alone.   Catching those two ideas, Romero et al. (1999) combined a genetic algorithm with a stochastic search procedure on icosahedrally derived lattices \cite{romerobg1999, barrongrs1999}.  The structures of the optimal  solutions  gotten  in \cite{romerobg1999} are either  icosahedral  or decahedral  in character. The  best  results  of \cite{wolf1998} for N  = 82, 84, 86, 88, 92, 93, 94, 95, 96, 99, 100 were gotten  by using a genetic  algorithm  alone.  Deaven  et al. (1996) also using the genetic  algorithm got the optimal value known for the magic number  $N=$ 88 \cite{deaventmh1996}. 
\item[$-$] The successful works to improve  Northby's results in \cite{northby1987} were mainly  done by Xue \cite{xue1994a,xue1994b}, Leary \cite{leary1997}, and Doye et al. \cite{doyewb1995,doyew1995}. 
\begin{enumerate}
\item[$\bullet$] Xue (1994a) introduced a modified version of the Northby  algorithm \cite{xue1994a}. He showed that in some cases the relaxation of the outer shell lattice local minimizer with a worse potential function value may lead to a local minimizer with a better value. In Northby's algorithm \cite{northby1987} the lattice  search  part is a  discrete  optimization  local  search  procedure, which  makes a lattice move to its neighboring  lattice with O($N^{\frac{5}{3}}$) time complexity. In \cite{xue1994a} Xue (1994a) introduced a simple storage  data structure to reduce the time complexity to O($N^{\frac{2}{3}}$) per move; and  then used a two-level  simulated  annealing  algorithm  within the supercomputer  CM-5 to be able to solve fastly the LJ problem  with sizes as large as 100,000 atoms.  In \cite{xue1994b} by employing  AVL trees \cite{horowitzs1990} data structure Xue (1994b) furthermore reduced  the time complexity  to O($\log N$ ) if NN (nearest  neighbor) potential function is used.  Xue (1994b) relaxed every lattice local minimizer found instead of relaxing only those lattice local minimizers with best known potential function value by a powerful Truncated Newton local search method \cite{xue1994b}, and at last got the best results known for $N=$ 65, 66, 134, 200, 300. 
\item[$\bullet$] Leary (1997) gave a successful Big Bang Algorithm \cite{leary1997} for getting the best values known of $N =$ 69, 78, 88, 107, 113, 115.  In \cite{leary1997}, the FCC  lattice  structure is discussed and  its  connections  are made  with  the macrocluster  problem. It is also concluded  in \cite{leary1997} that almost all known exceptions to global optimality of the well-known Northby multilayer icosahedral  conformations for microclusters are shown to be minor variants of that geometry.  The Big Bang Algorithm contains 3 steps:  Step 1 is an initial solution generating procedure  which randomly  generates each coordinate of the initial solution with the independently normal distribution; Step 2 is to generate the new neighborhood solution by discrete-typed  fixed step steepest descent method, which is repeated until no further progress is made;  Step 3 is to relax the best solution gotten  in Step 2 by a continuous optimization method--conjugate gradient method.
\item[$\bullet$] Doye et al. (1995) investigated the structures of clusters by mapping  the structure of the global minimum as a function of both cluster size and the range of the pair potential which is appropriate to the clusters of diatomic molecule, C$_{60}$  molecule, and  the ones between them both \cite{doyew1995}.  For the larger clusters the structure of the global minimum changes from icosahedral  to decahedral  to fcc as the range is decreased \cite{doyew1995}. In \cite{doyewb1995}, Doye et al. (1995) predicted the growth sequences for small decahedral  and fcc clusters by maximisation of the number  of NN contacts. 
\end{enumerate}
\item[$-$] Calvo et al. (2001) gave some results on quantum LJ Clusters in the use of Monte Carlo methods \cite{calvodw2001}.
\item[$-$] Xiang et al. (2004a) presented an efficient method based on lattice construction and the genetic algorithm and got global minima for $N =$ 310$\sim$561 \cite{xiangjcs2004}  In 2004, Xiang et al. (2004b) continued to present global minima for $N =$ 562$\sim$1000 \cite{xiangccs2004}.
\item[$-$] Barron-Romero (2005) found the best solutions for $N =$ 542--3, 546--8 in the use of a modified peeling greedy search method \cite{barronromero2005}.
\item[$-$] Takeuchi (2006) found best solutions for $N =$ 506, 521, 537--8 and 541 by a clever and efficient method ``using two operators: one modifies a cluster configuration by moving atoms to the most stable positions on the surface of a cluster and the other gives a perturbation on a cluster configuration by moving atoms near the center of mass of a cluster" \cite{takeuchi2006}.
\item[$-$] Lai et al. (2011a) found best solutions for $N =$ 533 and 536 using the dynamic lattice searching method with two-phase local search and interior operation \cite{laixh2011a, laixh2011b, yexh2011}.
\item[$-$] Algorithms to get the structures at the magic numbers $N=$ 17, 23, 24, 72, 88 (the exceptions to \cite{romerobg1999}):   
\begin{enumerate}
\item[$\bullet$] Freeman et al. (1985) presented the best value for $N=$ 17 when the thermodynamic properties of argon clusters were studied by a combination of classical and  quantum Monte Carlo  methods \cite{freemand1985}.  The  poly-icosahedral  growth of Farges  et al. (1985) starts from a 13-atom primitive icosahedron containing  a central  atom  and  12 surface atoms \cite{farges_etal1985}.  On each one of the five tetrahedral sites,  surrounding a particular  vertex,  a new atom  is added  and  finally a sixth  atom is placed  on top to create  a pentagonal  cap.   In this  way a 19-atom  structure being made  of double  interpenetrating icosahedra,  which is a 13-atom icosahedra  sharing  9 atoms, is obtained;  i.e., for three pentagonal bipyramids each one shares an apex with its nearest neighbour. In this way a 23-atom  model consisting of three interpenetrating  icosahedra  is gotten for the best value  known.   
\item[$\bullet$] Wille  (1987) used  the SA  method yielding low-lying energy states whose distribution depends  on the cooling rate to find the best solution known for $N=$ 24 \cite{wille1987}. 
\item[$\bullet$] Coleman  et al. (1997) proposed a build-up  process  to construct the optimal solution structures.  The  HOC  (half octahedral cap)  structure of the optimal solution for $N=$ 72 is found  by a prototype algorithm designed using the anisotropic effective  energy simualted annealing  method at each  build-up  stage (\cite{colemansw1997}).   
\item[$\bullet$] Wales  \& Doye  (1997) gave  the lowest values known  for $N =$ 192, 201 \cite{walesd1997}.  Their  method  is so-called  basin-hopping method,  in which first  the transformed  function  $\tilde{f} (x) = \min \{ f(x) \}$ was defined and  performed  starting from $x$ by the PR conjugate gradient method and then the energy  landscape for the function $\tilde{f} (x)$ was explored using a canonical  Monte  Carlo simulation.  
\item[$\bullet$] Leary (2000) has developed techniques for moving along sequences of local minima with decreasing energies to arrive at good candidates for global optima and got the best value known on $N=$ 185. 
\end{enumerate}
\end{enumerate}

\vskip 1cm

Now we have the outline of some successful optimization methods used to solve Eq.s (\ref{minimize_potential_energy})$\sim$(\ref{minimize_LJ_potential_energy}). We have found the hybrid idea of optimization methods was not emphasized very much (especially for solving Eq. (\ref{minimize_LJ_potential_energy})). Thus, in Section 2 of this paper we will emphasize the hybrid idea of optimization methods by introducing our own hybrid methods used to solve Eq.s (\ref{minimize_potential_energy})$\sim$(\ref{minimize_LJ_potential_energy}). Section 3 will present our recent results of applying the hybrid idea of SD and CG and SD again to do EM of some NMR and X-ray prion protein structures in the PDB Bank (www.rcsb.org); interesting findings will be reported in this Section. Why we choose prion proteins in this study is due to prions effect humans and almost all animals for a major public health concern (e.g. milks and meats we daily drink and eat). At last, in Section 4, we give a concluding remark on the effective and efficient hybrid idea of optimization methods.    

\section{The hybrid idea and some hybrid optimization methods}
In this Section, we use how to construct molecular structures of prion amyloid fibrils at AGAAAAGA segment as an example to illuminate the hybrid idea and some hybrid optimization methods we designed.

\vskip 1cm

Neurodegenerative amyloid diseases such as Alzheimer's, Parkinson's and Hungtington's all featured amyloid fibrils. Prions also cause a number of neurodegenerative diseases too. All these amyloid fibrils in 3-dimensional quaternary structure have 8 classes of steric zippers, with strong van der Waals interactions between $\beta$-sheets and hydrogen bonds between $\beta$-strands. Currently, there is no structural information about prion AGAAAAGA amyloid fibrils because of unstable, noncrystalline and insoluble nature of this region, though numerous laboratory experimental results have confirmed this region owning an amyloid fibril forming property (initially described in 1992 by Stanley B. Prusiner's group). We also did accurate computer calculations on this region and confirmed the amyloid fibril property in this palindrome region \cite{zhanghwwz2012, zhangz2013}.

\vskip 1cm

In \cite{zhang2011a},  we constructed three models, model 1 belongs to Class 7 (antiparallel, face=back, up-up) and models 2--3  belong to Class 1 (parallel, face-to-face, up-up) of steric zippers. The models were firstly optimized by SD and then followed by CG. SD has fast convergence but it is slow when close to minimums. CG is efficient but its gradient RMS and GMAX gradient do not have a nice convergence. When the models could not be optimized furthermore, we employed standard SA method (that simulates the annealing process of crystal materials with Monte Carlo property). After SA, we refined the models by SD and then CG again. SA is a global search optimization method \cite{zhang2014} that can make local optimal jump out of / escape from the local trap. We found the refinement results in a loss of potential energy nearly the same magnitude as that of SA; this implies to us SA is very necessary and very effective in our molecular modeling. Numerical results show to us the hybrid is very necessary, effective and efficient.

\vskip 1cm

When the gradient or its generalizations of the target/objective function $E(x)$ are very complex in form or they are not known, derivative-free methods benefit optimization problems. In \cite{zhangsw2011}, we introduced derivative-free discrete gradient (DG) method \cite{bagirovkm2014} into the derivative-free global search SA optimization method or genetic algorithms (GAs, which simulate the process natural competitive selection, crossover, and mutation of species), and designed hybrid methods SADG, GADG. In implementation, the hybrids of DG + SADG / GADG + DG were used, and at last SD+CG + SA + SD+CG of Amber package \cite{amber14} were used to refine the models. We found the hybrids work very well, and more precise best solutions for$N =$ 39, 40, 42, 48, 55, 75, 76, and 97 were found and their figures show that their structures are more stable than the ones currently best solutions known. We also found the hybrid of evolutionary computations with simulate annealing SA-SAES($\mu +\lambda $), SA-SACEP perform better than evolutionary computations or SA work alone \cite{zhang2011b}.

\vskip 1cm

Canonical dual theory in some sense is the hybrid of the primal and the dual. In \cite{zhanggy2011}, we solved the dual problem and then got the solutions for the primal problem. We found the refinement using AMBER package is not necessary. This implies to us the hybrid of primal and dual in canonical dual theory is good enough and effective.

\vskip 1cm

As said in Section 1, in some cases, CG cannot be used to find the EM path; this point will also be shown in next Section (see Tab. \ref{energy_variations_during_energy_minimizations}). Thus, in \cite{zhanghwwz2012}, we specially studied and implement the LBFGS method designed by us and then hybridize it with the LBFGS method of AMBER package. We found the hybrid is very necessary and effective.

\vskip 1cm

By our numerical experiences shown above, the hybrid idea is very necessary, effective and efficient for some hybrid optimization methods in known packages or designed by us. In next Section, we will apply the hybrid idea to do some practical works for some important prion protein NMR and X-ray structures deposited in the PDB Bank. 

\section{An application to prion protein structures, focusing on the $\beta$2-$\alpha$2 loop}
Before we use the structure taken from PDB Bank, usually we need to relax it, in order to remove bad contacts and also fix up hydrogen positions. Fairly short  local optimization is sufficient to refine and relax the structure. We will use SD-CG-SD local optimization methods. In SD, its search direction is a n-dimensional search and its step-length search is a 1-dimensional search. In CG, the search is usually in a 2-dimensional subspace and conjugacy is a good property only associated with exact line search \cite{sunz2001}. Using the hybrid of SD and CG is also in order to remove all these (dimensional) unbalances. In our EM here, the free package Swiss-PdbViewer 4.1.0 (spdbv.vital-it.ch) that has been developed for 20 years is used, we set 3000 steps for SD, then 3000 steps for CG, and then 3000 steps for SD again, Bonds, Angles, Torsions, Improper, Non-bonded and Electrostatic are considered, 12.000 $\mathring{\text{A}}$ is chosen for the Cutoff, stop SD or CG when delta E between two steps is below 0.005 kJ/mol, and stop SD or CG when Force acting on any atom is below 1.000.

\vskip 1cm

We found, for the research of prion proteins, the S2-H2 loop (and its interactions with the C-terminal of H3) is a focus \cite{biljan_etal2012a,biljan_etal2012b,biljanigrzpl2011,calzolai_etal2000, christen_etal2009, christen_etal2012, dambergercphw2011, gossertblfw2005, ile_etal2010, leeahkssy2010, wenlxpyhl2010, wenlyxpxl2010, kong_etal2013, perezdw2010, perezw2008, sweeting_etal2013, zahngvsw2003, zhangszss2000, kurt_etal2014a, kurt_etal2014b, huangc2015}. All prion protein structures have high similarity in three $\alpha$-helices (H1, H2, H3) and two $\beta$-strands (S1, S2), but there is a great difference just at this S2-H2 loop:
\begin{enumerate}
\item[] (i) structure with disordered S2-H2 loop:
\begin{enumerate}
\item[$\bullet$] mousePrP (1AG2.pdb at 25 $\mathring{\text{}}$C), 
\item[$\bullet$] humanPrP (1QLX.pdb), 
\item[$\bullet$] bovinePrP (1DWY.pdb), 
\item[$\bullet$] SyrianHamsterPrP (1B10.pdb), 
\item[$\bullet$] dogPrP (1XYK.pdb) (- resist to prion infection), 
\item[$\bullet$] catPrP (1XYJ.pdb), 
\item[$\bullet$] sheepPrP (1UW3.pdb), 
\item[$\bullet$] mousePrP[N174T] (1Y15.pdb), 
\item[$\bullet$] humanPrP[Q212P]-M129 (2KUN.pdb), 
\item[$\bullet$] humanPrP-M129 (1QM1.pdb), 
\item[$\bullet$] rabbitPrP[S173N] (2JOH.pdb), 
\item[$\bullet$] rabbitPrP[I214V] (2JOM.pdb), 
\item[$\bullet$] rabbitPrP[S170N] (4HLS.pdb), 
\item[$\bullet$] rabbitPrP[S174N] (4HMM.pdb), 
\item[$\bullet$] rabbitPrP[S170N,S174N] (4HMR.pdb), 
\end{enumerate}
\item[] (ii) structure with highly and clearly ordered S2-H2 loop:
\begin{enumerate}
\item[$\bullet$] mousePrP (2L39.pdb at 37 $\mathring{\text{}}$C), 
\item[$\bullet$] mousePrP[V166A] (2KFO.pdb), 
\item[$\bullet$] mousePrP[D167S] (2KU5.pdb at 20 $\mathring{\text{}}$C),  
\item[$\bullet$] mousePrP[D167S,N173K] (2KU6.pdb), 
\item[$\bullet$] mousePrP[Y169G] (2L1D.pdb), 
\item[$\bullet$] mousePrP[Y169A] (2L40.pdb), 
\item[$\bullet$] mousePrP[S170N] (2K5O.pdb), 
\item[$\bullet$] mousePrP[S170N,N174T] (1Y16.pdb), 
\item[$\bullet$] mousePrP[F175A] (2L1E.pdb), 
\item[$\bullet$] mousePrP[Y225A,Y226A] (2KFM.pdb), 
\item[$\bullet$] mousePrP[Y169A,Y225A,Y226A] (2L1K.pdb at 20 $\mathring{\text{}}$C), 
\item[$\bullet$] elkPrP (1XYW.pdb), 
\item[$\bullet$] pigPrP (1XYQ.pdb), 
\item[$\bullet$] BankVolePrP (2K56.pdb), 
\item[$\bullet$] TammarWallabyPrP (2KFL.pdb), 
\item[$\bullet$] rabbitPrP (2FJ3.pdb, 3O79.pdb), 
\item[$\bullet$] horsePrP (2KU4.pdb),
\end{enumerate}
\end{enumerate}
where elk and Bank Vole can be infected by prions though they have a highly and clearly ordered S2-H2 loop, and the codes in the brackets are the PDB codes in the PDB Bank. For all these NMR and X-ray structures we did SD-CG-SD relaxation and the variations of the EMs are listed in Tab. \ref{energy_variations_during_energy_minimizations}. From Tab. \ref{energy_variations_during_energy_minimizations}, we can see the energy decreases from SD to CG and from CG to SD. For mousePrP[Y169A], CG is not working well, but it adjusts the SD methods so that it make SD-CG-SD work very well in the second round.  
\begin{table}[h!]
\caption{\textsf{Energy variations during the energy minimizations (number of iterations are in the brackets):}}
\centering
{\tiny
\begin{tabular}{|l                             |c                    |c                    |c|} \hline
                Species                       &SD                   &CG                   &SD\\ \hline \hline
                mousePrP                      &-7234.507 (391)      &-7243.597 (13)       &-7275.160 (45)\\  \hline
                humanPrP                      &-7460.885 (296)      &-7610.384 (195)      &-7640.157 (54)\\  \hline
                bovinePrP                     &-7698.001 (544)      &-7809.687 (255)      &-7819.315 (26)\\  \hline
                SyrianHamsterPrP              &-7418.702 (225)      &-7653.251 (258)      &-7688.044 (50)\\  \hline
                dogPrP                        &-7148.517 (549)      &-7225.084 (151)      &-7251.191 (60)\\  \hline
                catPrP                        &-6935.915 (192)      &-7186.961 (186)      &-7382.361 (229)\\  \hline
                sheepPrP                      &-8066.183 (291)      &-8066.300 (22)       &-8204.010 (179)\\  \hline
                mousePrP[N174T]               &-7418.443 (211)      &-7657.886 (206)      &-7923.308 (419)\\  \hline
                humanPrP[Q212P]-M129          &-7662.032 (464)      &-7688.836 (49)       &-7846.123 (424)\\  \hline
                humanPrP-M129                 &-7249.798 (350)      &-7345.074 (114)      &-7410.470 (99)\\  \hline
                rabbitPrP[S173N]-NMR          &-7173.492 (271)      &-7666.456 (492)      &-7684.357 (36)\\  \hline
                rabbitPrP[I214V]-NMR          &-7785.640 (774)      &-7802.710 (39)       &-7835.354 (88)\\  \hline
                rabbitPrP[S170N]-X-ray        &-8682.104 (414)      &-8753.827 (178)      &-8827.220 (254)\\  \hline
                rabbitPrP[S174N]-X-ray        &-8551.921 (286)      &-8659.515 (160)      &-8739.475 (178)\\  \hline
                rabbitPrP[S170N,S174N]-X-ray  &-8864.615 (363)      &-8899.511 (65)       &-8945.262 (106)\\ \hline \hline
                mousePrP - at 37 $\mathring{\text{}}$C             &-7679.395 (212)      &-8031.103 (397)      &-8134.173 (213)\\  \hline
                mousePrP[V166A]               &-8040.172 (436)      &-8153.529 (174)      &-8164.996 (22)\\  \hline
                mousePrP[D167S] - at 20 $\mathring{\text{}}$C      &-7938.545 (594)      &-7987.788 (83)       &-8051.072 (121)\\  \hline
                mousePrP[D167S,N173K]         &-7615.915 (546)      &-7751.804 (205)      &-7906.153 (290)\\  \hline
                mousePrP[Y169G]               &-7713.983 (249)      &-7885.343 (147)      &-7913.277 (27)\\  \hline
                mousePrP[Y169A]               &-7948.267 (507)      &-7949.964 (1)        &-7955.516 (3)\\ 
                                              &-7972.101 (17)       &-8064.175 (142)      &-8143.807 (179)\\ \hline 
                mousePrP[S170N]               &-7341.668 (74)       &-7790.935 (262)      &-7947.893 (224)\\  \hline
                mousePrP[S170N,N174T]         &-7988.545 (468)      &-8126.978 (241)      &-8149.948 (55)\\  \hline
                mousePrP[F175A]               &-7660.438 (680)      &-7859.632 (381)      &-7873.355 (34)\\  \hline
                mousePrP[Y225A,Y226A]         &-7457.356 (244)      &-7588.338 (109)      &-7757.797 (213)\\  \hline
                mousePrP[Y169A,Y225A,Y226A] - at 20 $\mathring{\text{}}$C &-7609.773 (173) &-7690.231 (45)     &-7700.162 (4)\\  \hline
                elkPrP                        &-7894.305 (875)      &-7959.371 (160)      &-7978.686 (53)\\  \hline
                pigPrP                        &-6354.886 (141)      &-6735.155 (321)      &-6813.048 (109)\\  \hline
                bankVolePrP                   &-7727.260 (478)      &-7799.991 (118)      &-7951.502 (376)\\  \hline
                tammarWallabyPrP              &-8028.082 (393)      &-8195.817 (248)      &-8238.248 (69)\\  \hline
                rabbitPrP-NMR                 &-7712.972 (814)      &-7730.035 (36)       &-7790.620 (154)\\  \hline
                rabbitPrP-X-ray               &-8643.309 (602)      &-8658.327 (40)       &-8717.344 (187)\\  \hline
                horsePrP                      &-7335.273 (226)      &-7614.526 (433)      &-7633.569 (45)\\  \hline     
\end{tabular} 
} \label{energy_variations_during_energy_minimizations}
\end{table}

\begin{table}[h!]
\caption{\textsf{Hydrogen bonds at the S2-H2 loop, and its linkage with the C-terminal end of H3 (in the brackets are the distances of the hydrogen bonds):}}
\centering
{\tiny
\begin{tabular}{|l                            |l                                   |l|} \hline
                Species                       &at the S2-H2 loop                   &linking with the C-terminal end of H3\\ \hline \hline
                mousePrP                      &D178.OD1--R164.2HH1 (1.97 $\mathring{\text{A}}$) &S170:OG--Y225:HH (1.71 $\mathring{\text{A}}$)\\ \hline
                humanPrP                      & &\\ \hline
                bovinePrP                     & &\\ \hline
                SyrianHamsterPrP              &Y169.HH--D178.OD2 (1.80 $\mathring{\text{A}}$) &\\
                                              &T183.OG1--Y162.H (1.74 $\mathring{\text{A}}$) &\\ \hline
                dogPrP                        &D167.O--S170.H (1.82 $\mathring{\text{A}}$) &\\ \hline
                catPrP                        & &\\ \hline
                sheepPrP                      &P165.O--Q168.H (1.90 $\mathring{\text{A}}$)     &Y163.OH--Q217.2HE2 (1.91 $\mathring{\text{A}}$)\\
                                              &D178.OD2--Y128.HH (1.81 $\mathring{\text{A}}$)  &\\
                                              &Y163.O--M129.H (1.79 $\mathring{\text{A}}$) &\\
                                              &N171.O--N171.2HD2 (1.84 $\mathring{\text{A}}$)  &\\
                                              &T183.OG1--Y162.H (1.92 $\mathring{\text{A}}$)  &\\ \hline
                mousePrP[N174T]               &D178.OD2--Y128.HH (1.72 $\mathring{\text{A}}$)  &\\
                                              &Y163.H--M129.O (1.79 $\mathring{\text{A}}$)     &\\ \hline
                humanPrP[Q212P]-M129          &N171.OD1--N173.H (1.84 $\mathring{\text{A}}$)   &Q172.OE1--Y218.HH (1.68 $\mathring{\text{A}}$)\\ 
                                              &Y162.H--T183.OG1 (1.94 $\mathring{\text{A}}$)   &\\ \hline
                humanPrP-M129                 &Y162.H--T183.OG1 (1.89 $\mathring{\text{A}}$)   &\\ \hline
                rabbitPrP[S173N]-NMR          &                                                    &E220.OE1--Y162.HH (1.68 $\mathring{\text{A}}$)\\ 
                                              &                                                    &D177.OD2--Y127.HH (1.71 $\mathring{\text{A}}$)\\ \hline
                rabbitPrP[I214V]-NMR          &Y168.O--N170.2HD2 (1.86 $\mathring{\text{A}}$)  &\\
                                              &Q171.O--V175.H (1.70 $\mathring{\text{A}}$)     &\\ \hline
                rabbitPrP[S170N]-X-ray        &P165.O--Q168.H (1.93 $\mathring{\text{A}}$)     &Q172.OE1--Q219.2HE2 (1.85 $\mathring{\text{A}}$)\\
                                              &V166.O--Y169.H (1.92 $\mathring{\text{A}}$)     &\\ \hline
                rabbitPrP[S174N]-X-ray        &P165.O--Q168.H (1.94  $\mathring{\text{A}}$)    &\\ 
                                              &V166.O--Y169.H (1.92 $\mathring{\text{A}}$)     &\\
                                              &Y169.HH--D178.OD2 (1.74 $\mathring{\text{A}}$)  &\\
                                              &N171.HD2--N174.ND2 (1.95 $\mathring{\text{A}}$) &\\ \hline
                rabbitPrP[S170N,S174N]-X-ray  &P165.O--Q168.H (1.92 $\mathring{\text{A}}$)     &Q172.OE1--Q219.2HE2 (1.89 $\mathring{\text{A}}$)\\ 
                                              &V166.O--Y169.H (1.93 $\mathring{\text{A}}$)     &\\
                                              &Y169.HH--D178.OD2 (1.78 $\mathring{\text{A}}$)   &\\
                                              &N171.OD1--N174.H (1.88 $\mathring{\text{A}}$)   &\\ \hline \hline
                mousePrP - at 37 $\mathring{\text{}}$C &Y163.O--M129.H (1.71 $\mathring{\text{A}}$)  &S170.O--Y218.HH (1.61 $\mathring{\text{A}}$)\\
                                                       &R164.O--Y169.HH (1.80 $\mathring{\text{A}}$) &\\
                                                       &V166.O--Y169.H (1.83 $\mathring{\text{A}}$)  &\\
                                                       &D178.OD2--Y128.HH (1.74 $\mathring{\text{A}}$) &\\ \hline
                mousePrP[V166A]               &R164.HE--Q168.OE1 (1.90 $\mathring{\text{A}}$) &\\
                                              &T183.OG1--Y162.H (2.02  $\mathring{\text{A}}$) &\\ \hline
                mousePrP[D167S] - at 20 $\mathring{\text{}}$C      &R164.2HH2--Q168.OE1 (1.90 $\mathring{\text{A}}$) &\\
                                                                   &V166.O--Y169.H (1.87 $\mathring{\text{A}}$)      &\\
                                                                   &T183.OG1--Y162.H (1.90 $\mathring{\text{A}}$)    &\\  \hline
                mousePrP[D167S,N173K]         &R164.CZ--Q168.1HE2 (1.96 $\mathring{\text{A}}$) &\\
                                              &V166.O--Y169.H (1.82 $\mathring{\text{A}}$)     &\\
                                              &V166.O--S170.H (1.75 $\mathring{\text{A}}$)     &\\ \hline
                mousePrP[Y169G]               &T183.OG1--Y162.H (1.95 $\mathring{\text{A}}$) &\\ \hline
                mousePrP[Y169A]               &P165.O--S170.HG (1.91 $\mathring{\text{A}}$) &Q168.O--Y225.HH (1.73 $\mathring{\text{A}}$)\\
                                              &T183.OG1--Y162.H (1.94 $\mathring{\text{A}}$) &\\ \hline
                mousePrP[S170N]               &N170.2HD2--N171.OD1 (1.91 $\mathring{\text{A}}$) &\\
                                              &T183.OG1--Y162.H (1.84 $\mathring{\text{A}}$) &\\ \hline
                mousePrP[S170N,N174T]         &Y169.HH--D178.OD2 (2.01 $\mathring{\text{A}}$) &\\
                                              &N171.O--N171.2HD2 (1.81 $\mathring{\text{A}}$) &\\
                                              &T183.OG1--Y162.H (1.79  $\mathring{\text{A}}$) &\\ \hline
                mousePrP[F175A]               &V166.O--Y169.H (1.76 $\mathring{\text{A}}$)      &R164.O--Y218.HH (1.61 $\mathring{\text{A}}$)\\
                                              &N171.1HD2--N174.ND2 (1.75 $\mathring{\text{A}}$) &\\
                                              &N171.2HD2--N171.O (1.78 $\mathring{\text{A}}$)   &\\ \hline
                mousePrP[Y225A,Y226A]         &D178.OD2--Y128.HH (1.61 $\mathring{\text{A}}$)   &\\
                                              &Y163.O--MET129.H (1.77 $\mathring{\text{A}}$)      &\\
                                              &P165.O--Q168.2HE2 (1.85 $\mathring{\text{A}}$)   &\\
                                              &S170.OG--N174.1HD2 (1.95 $\mathring{\text{A}}$)  &\\
                                              &N171.H--N174.OD1 (1.94 $\mathring{\text{A}}$)    &\\
                                              &N173.2HD2--HIS177.NE2 (1.93 $\mathring{\text{A}}$) &\\
                                              &T183.OG1--Y162.H (1.96 $\mathring{\text{A}}$) &\\ \hline
                mousePrP[Y169A,Y225A,Y226A] - at 20 $\mathring{\text{}}$C & &\\ \hline
                elkPrP                        &P165--Q168.H (1.80 $\mathring{\text{A}}$) &\\
                                              &Y169.HH--D178.OD2 (1.73 $\mathring{\text{A}}$) &\\
                                              &T183.OG1--Y162.H (1.79 $\mathring{\text{A}}$) &\\ \hline
                pigPrP                        & &\\ \hline
                BankVolePrP                   &P165.O--Q168.H (1.99 $\mathring{\text{A}}$) &\\
                                              &T183.OG1--Y162.H (1.92 $\mathring{\text{A}}$) &\\ \hline
                TammarWallabyPrP              &P165.O--Q168.N (1.98 $\mathring{\text{A}}$) &\\
                                              &I166.O--Y169.H (1.92 $\mathring{\text{A}}$) &\\
                                              &T183.OG1--M162.H (1.79 $\mathring{\text{A}}$) &\\ \hline
                rabbitPrP-NMR                 &T182.OG1--Y161.H (1.81 $\mathring{\text{A}}$) &\\ \hline
                rabbitPrP-X-ray               &V166.O--Y169.H (1.95 $\mathring{\text{A}}$) &\\
                                              &D178.OD2--Y169.HH (1.82 $\mathring{\text{A}}$) &\\ \hline
                horsePrP                      &D178.OD2--Y169.HH (1.71 $\mathring{\text{A}}$) &\\ \hline      
\end{tabular} 
} \label{hydorgen_bonds_in_loop_optimized}
\end{table}   

\begin{table}[h!]
\caption{\textsf{Salt bridges at the S2-H2 loop, and its linkage with the C-terminal end of H3 (in the brackets are the distances of the salt bridges):}}
\centering
{\tiny
\begin{tabular}{|l                            |l                                   |l|} \hline
                Species                       &at the S2-H2 loop                   &linking to the C-terminal end of H3\\ \hline \hline
                mousePrP                      &D178.OD1--R164.NE (4.01 $\mathring{\text{A}}$)   &\\
                                              &D178.OD1--R164.NH1 (2.96 $\mathring{\text{A}}$)  &\\ \hline
                humanPrP                      & &\\ \hline
                bovinePrP                     & &\\ \hline
                SyrianHamsterPrP              & &\\ \hline
                dogPrP                        & &\\ \hline
                catPrP                        &D178.OD1--R164.NH2 (4.59 $\mathring{\text{A}}$)  &\\ 
                                              &D178.OD2--R164.NE (2.97 $\mathring{\text{A}}$)   &\\
                                              &D178.OD2--R164.NH1 (4.76 $\mathring{\text{A}}$)  &\\
                                              &D178.OD2--R164.NH2 (3.06 $\mathring{\text{A}}$)  &\\ \hline
                sheepPrP                      & &\\ \hline
                mousePrP[N174T]               & &\\ \hline
                humanPrP[Q212P]-M129          &E168.OE1--R164.NH1 (3.07 $\mathring{\text{A}}$)  &\\ 
                                              &E168.OE1--R164.NH2 (3.33 $\mathring{\text{A}}$)  &\\ \hline
                humanPrP-M129                 &D178.OD2--R164.NH1 (2.99 $\mathring{\text{A}}$)  &\\ 
                                              &D178.OD2--R164.NH2 (3.13 $\mathring{\text{A}}$)  &\\ 
                                              &D178.OD2--R164.NE (4.75 $\mathring{\text{A}}$)   &\\
                                              &D178.OD1--R164.NH1 (4.61 $\mathring{\text{A}}$)  &\\ \hline
                rabbitPrP[S173N]-NMR          & &\\ \hline
                rabbitPrP[I214V]-NMR          & &\\ \hline
                rabbitPrP[S170N]-X-ray        &D178.OD2--R164.NH2 (3.16 $\mathring{\text{A}}$)  &\\ 
                                              &D178.OD2--R164.NE (4.14 $\mathring{\text{A}}$)   &\\
                                              &D178.OD1--R164.NH2 (4.06 $\mathring{\text{A}}$)  &\\ 
                                              &D178.OD1--H177.ND1 (2.95 $\mathring{\text{A}}$)  &\\ \hline
                rabbitPrP[S174N]-X-ray        & &\\ \hline
                rabbitPrP[S170N,S174N]-X-ray  & &\\ \hline \hline
                mousePrP - at 37 $\mathring{\text{}}$C             & &\\ \hline
                mousePrP[V166A]               & &\\ \hline
                mousePrP[D167S] - at 20 $\mathring{\text{}}$C      & &\\ \hline
                mousePrP[D167S,N173K]         & &\\ \hline
                mousePrP[Y169G]               & &\\ \hline
                mousePrP[Y169A]               & &\\ \hline
                mousePrP[S170N]               & &D167.OD1--R229.NH1 (3.03 $\mathring{\text{A}}$)\\
                                              & &D167.OD1--R229.NH2 (3.81 $\mathring{\text{A}}$)\\ \hline
                mousePrP[S170N,N174T]         & &\\ \hline
                mousePrP[F175A]               & &\\ \hline
                mousePrP[Y225A,Y226A]         & &\\ \hline
                mousePrP[Y169A,Y225A,Y226A] - at 20 $\mathring{\text{}}$C & &\\ \hline
                elkPrP                        & &\\ \hline
                pigPrP                        &D178.OD2--R164.NH1 (3.00 $\mathring{\text{A}}$) &\\ \hline
                BankVolePrP                   & &\\ \hline
                TammarWallabyPrP              & &\\ \hline
                rabbitPrP-NMR                 & &H176.ND1--E210.OE1 (2.96 $\mathring{\text{A}}$)\\ \hline
                rabbitPrP-X-ray               & &\\ \hline
                horsePrP                      & &\\ \hline      
\end{tabular} 
} \label{salt_bridges_in_loop_optimized}
\end{table}   

\begin{figure}[h!] 
\centerline{
\includegraphics[width=5.2in]{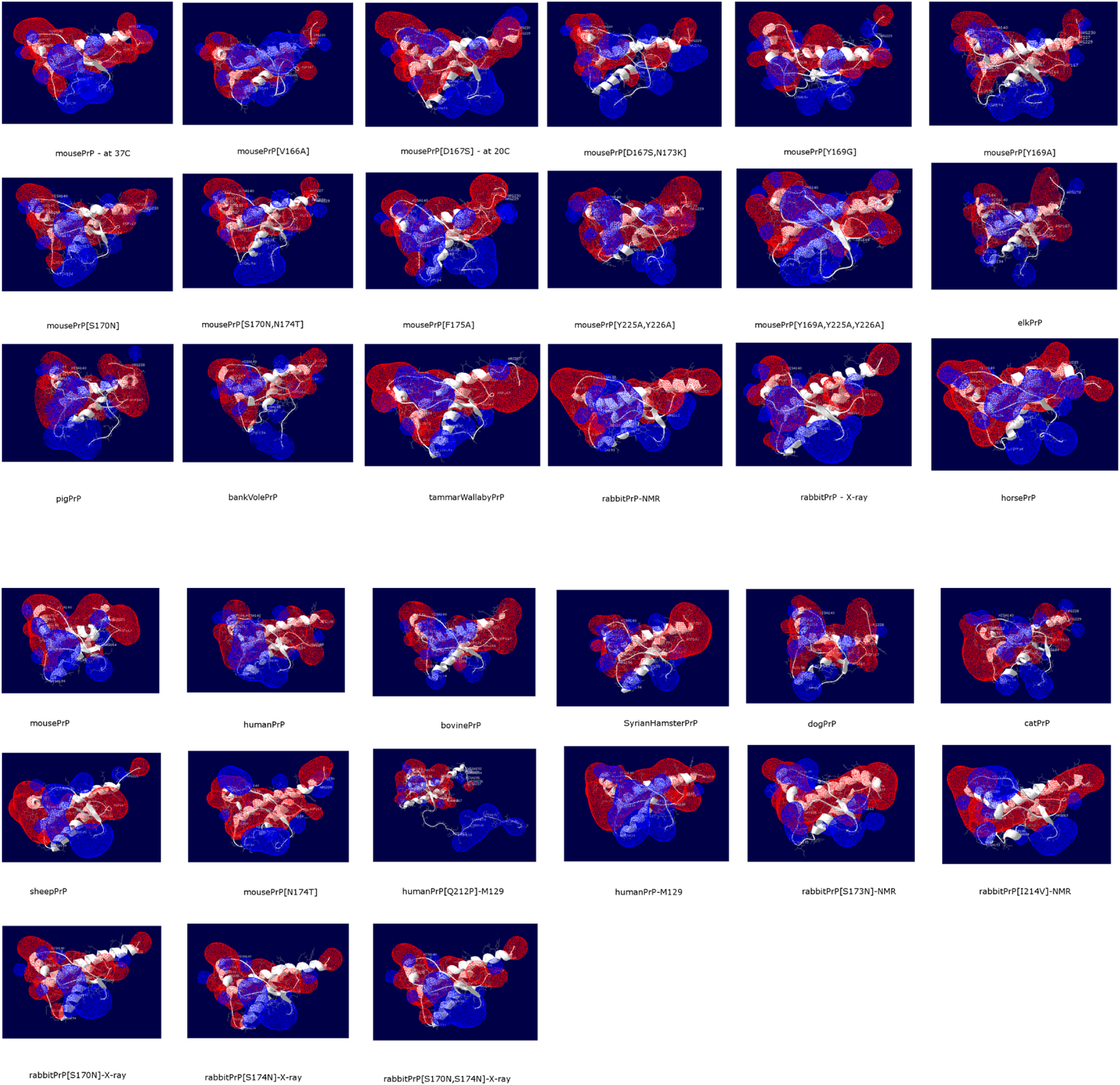}
}
\caption{\textsf{Positively (in blue) and negatively (in red) charged residues distributed on each protein structure surface of the 33 PrPs. The first three rows of PrPs are with highly and clearly ordered S2--H2 loop; but the last three rows of PrPs are with disordered S2--H2 loop.}} \label{charge_distributions}
\end{figure}
\begin{figure}[h!] 
\centerline{
\includegraphics[width=5.2in]{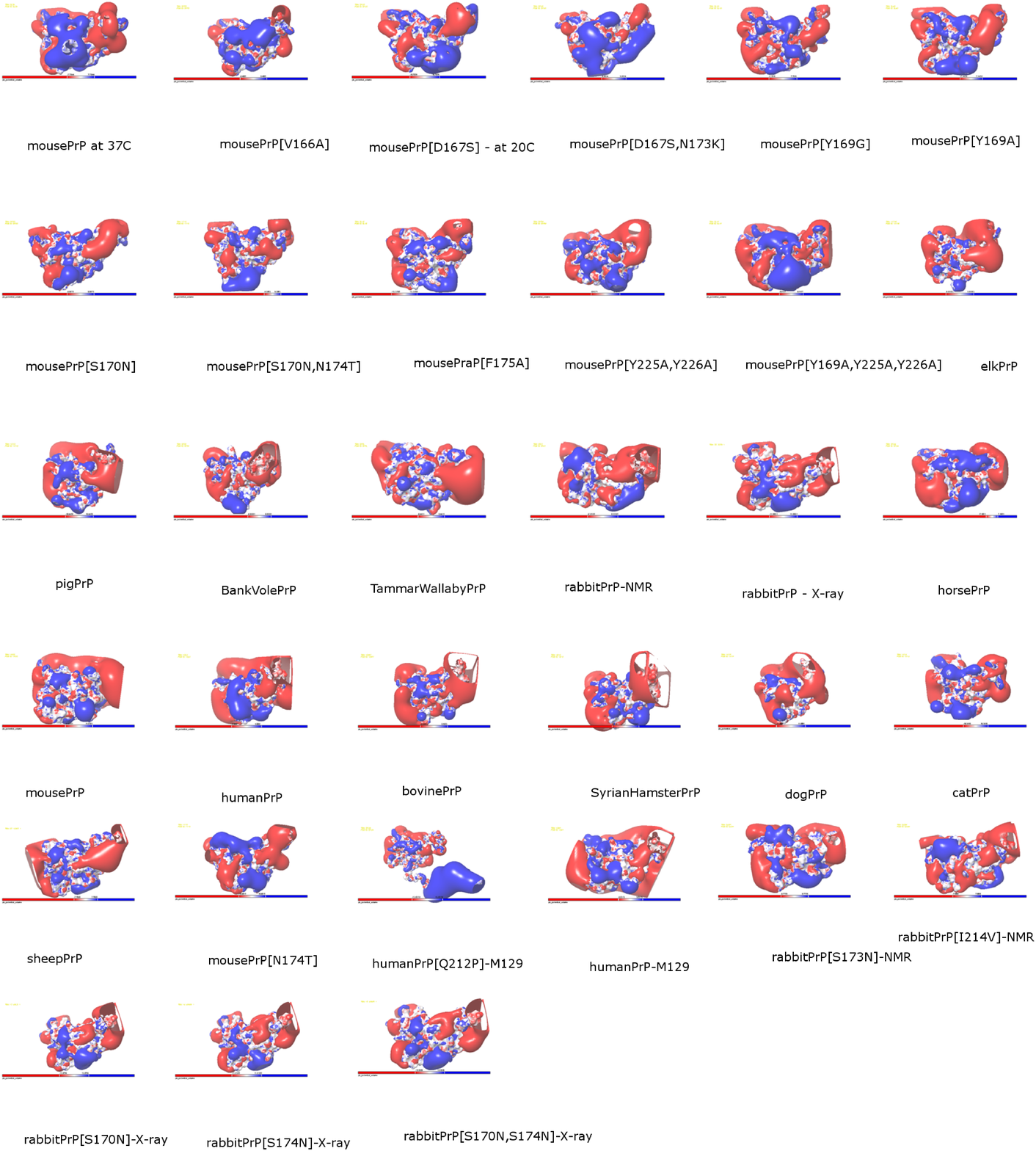}
}
\caption{\textsf{Surface electrostatic charge distributions for each of the 33 PrPs. The first three rows of PrPs are with highly and clearly ordered S2--H2 loop; but the last three rows of PrPs are with disordered S2--H2 loop. Blue is for positive charge whereas red is for negative charge.}} \label{surface_electrostatic_distributions}
\end{figure}
After the SD-CG-SD relaxation of all the structures, now these optimized structures can be used to obtain some helpful structural information [e.g. (i) hydrogen bonds (see Tab. \ref{hydorgen_bonds_in_loop_optimized}), (ii) electrostatic charge distributions on the protein structure surface (see Fig. \ref{charge_distributions} and Tab. \ref{cliques_of_positively_charged_residues}),  (iii) salt bridges (see Tab. \ref{salt_bridges_in_loop_optimized}), and (iv) $\pi$-$\pi$-stacking and $\pi$-cations (see Tab. \ref{pi_stacks}); here why we consider the information of (i)$\sim$(iv) is due to ``the performance of protein biological function is driven by a number of non-covalent interactions such as hydrogen bonding, ionic interactions, Van der Waals forces, and hydrophobic packing" (en.wikipedia.org/wiki/Protein\_structure)] at the S2-H2 loop, in order to furthermore understand the S2-H2 loop: why some species has a clearly and highly ordered S2-H2 loop and why some species just has a disordered S2-H2 loop. (i) Using the VMD package (www.ks.uiuc.edu/Research/vmd/), with Tab. \ref{hydorgen_bonds_in_loop_optimized}, we may observed that the species owning the disordered S2-H2 loop usually does not have a  3$_{10}$-helix in the S2-H2 loop (except for sheepPrP, rabbitPrP[S170N]-X-ray, rabbitPrP[S174N]-X-ray, and rabbitPrP[S170N,S174N]-X-ray), but the species that has the clearly and highly ordered S2-H2 loop usually owns a 3$_{10}$-helix, constructed by the following hydrogen bond(s) respectively:
\begin{enumerate}
\item[$\bullet$] V166--Y169 - mousePrP at 37 $\mathring{\text{}}$C, mousePrP[D167S,N173K], mousePrP[F175A], rabbitPrP-X-ray,
\item[$\bullet$] R164--Q168 - mousePrP[V166A], 
\item[$\bullet$] P165--Q168 - mousePrP[Y225A,Y226A], elkPrP, BankVolePrP, 
\item[$\bullet$] P165--Q168, I166--Y169 - TammarWallabyPrP.
\end{enumerate}
(ii) Seeing Fig. \ref{charge_distributions} and Tab. \ref{cliques_of_positively_charged_residues}, we may know that at the S2-H2 loop it is mainly covered by the electrical cloud of negatively (in red color) charged residues [except for mousePrP[D167S,N173K], rabbitPrP-NMR and horsePrP etc (Fig. \ref{surface_electrostatic_distributions})], with positively (in blue) charged residues R164 (for all species) and H177 (for mousePrP, humanPrP-M129, mousePrP[D167S,N173K] only) at the N-terminal end and C-terminal end of S2-H2 loop respectively (we found there is a salt bridge R164--D178 linking this loop of rabbitPrP-NMR, rabbitPrP-X-ray, horsePrP, dogPrP, elkPrP and buffaloPrP for long time MD simulations \cite{zhang2010, zhang2011d, zhang2011e, zhangl2011, zhangz2014, zhangwz2015, zhang2015a, zhang2011f}). From Tab. \ref{cliques_of_positively_charged_residues}, we might see that the negatively charged S2-H2 loop might have long distance nuclear overhauser effect (NOE) interactions with the positively charged residues such as K204, R208, K/R220, R227, R228, R229, and R230 at the C-terminal end of H3. (iii) The salt bridges in Tab. \ref{salt_bridges_in_loop_optimized} might be not very strong and will be quickly broken in a long time MD simulations \cite{zhang2011e}. (iv) Lastly, we present some bioinformatics of $\pi$-$\pi$-stacking and $\pi$-cations (one kind of van der Waals interactions) at the S2-H2 loop. Seeing Tab. \ref{pi_stacks}, we may know at S2-H2 loop and its contacts with the C-terminal end of H3 there are the following $\pi$-$\pi$-stacks Y169--F175, F175--Y218, Y163--Y218, and the following $\pi$-cations R164--Y169, R164--Y128, which clearly contribute to the clearly and highly ordered S2-H2 loop structures \cite{zhang2015b}. For buffaloPrP, we found another two $\pi$-stackings: Y163--F175--Y128 \cite{zhang2015b, zhang2015a}. Thus, for PrPs, we found an interesting ``$\pi$-circle" Y128--F175--Y218--Y163--F175--Y169--R164--Y128(--Y162) around the S2-H2 loop.
\begin{table}[h!]
\caption{\textsf{Cliques of positively charged residues distributed on the surface of each optimized structure of the 33 PrPs (one bracket is one clique):}}
\centering
{\tiny
\begin{tabular}{l                            |l} \hline
mousePrP & (R164), (H177), (K220), (R156), (R148), (K185), (R151, R136, K194)\\ \hline
humanPrP & (R164), (R228), (R220), (H155, R156, K185, H187, K194,)\\ \hline
bovinePrP & (R164), (K185), (K194), (H155, R156, R136, H187, R220)\\ \hline
SyrianHamsterPrP & (R164), (K185), (K194), (K220), (R136, R156), (R148, R151)\\ \hline
dogPrP & (R164), (R228), (R220), (K185), (K194), (R151, R156)\\ \hline
catPrP & (R164), (H187, R185), (K194), (R228), (R229), (K220)\\ \hline
sheepPrP & (R164), (K185), (R228), (H187, K194, H155, R156, R136, R220)\\ \hline
mousePrP[N174T] & (R164), (R229, R230), (K220), (K185, H187, K194) (R156), (H140, R151, R136, R148)\\ \hline
humanPrP[Q212P]-M129 & (R164), (K185), (R228)\\ \hline
humanPrP-M129 & (R164, H177, K185, H187, K194, H155, R156, R220, R208, R151, R148, H140), (R228)\\ \hline
rabbitPrP[S173N]-NMR & (R163, R135, K184, H186), (R227), (K193, R155), (R147), (H139, R135, R150)\\ \hline
rabbitPrP[I214V] & (R163, R135), (R227), (H186, K184), (K193), (R147), (H139, R135, R150)\\ \hline
rabbitPrP[S170N]-X-ray & (R164, R136, H187, K185, K194), (R228), (H140, R208), (H140), (R148), (R156)\\ \hline
rabbitPrP[S174N]-X-ray & (R164, R136, K185, - H187, R156, K194, R136), (R228), (H140), (R148), (R156)\\ \hline
rabbitPrP[S170N,S174N]-X-ray & (R164, R136, K185, H187, K194, R156, R136), (R228), (H140), (R151), (R148), (R156)\\ \hline \hline
mousePrP at 37 $\mathring{\text{}}$C & (R164), (R229, R230), (K220), (K185, H187, K194, R136), (R148), (R151)\\ \hline
mousePrP[V166A] & (R164, K185), (R136, K220, R229, R208), (R230), (K194), (R156), (R148), (R151)\\ \hline
mousePrP[D167S] at 20 $\mathring{\text{}}$C & (R164, K185, H187, K194), (R229, R230), (R156), (R151, R136), (R148), (K220)\\ \hline
mousePrP[D167S,N173K] & (R164), (H177), (R229 R230), (K185, H187), (K194), (R136, R156), (R156), (R148),(R151),(K220)\\ \hline
mousePrP[Y169G] & (R164, K185, K194), (R230), (R229), (K220), (R136), (R156), (R151), (R148)\\ \hline
mousePrP[Y169A] & (R164, K185, K194), (R230), (R229), (K220), (R136, R151), (R156), (R148)\\ \hline
mousePrP[S170N] & (R164), (R230), (R229), (K185), (K194, H187, R136), (R156), (K220), (R208), (R151), (R148)\\ \hline
mousePrP[S170N,N174T] & (R164), (R229), (R230), (K194, K185), (R156), (R136), (R148), (R151), (R208), (K220)\\ \hline
mousePrP[F175A] & (R164, K185), (R229), (R230), (K194, H187, R136, R156, R151), (R148), (H140)\\ \hline
mousePrP[Y225A,Y226A] & (R164, K185), (R229, R230), (R156), (K194, R136), (R148, R151), (R208), (K220)\\ \hline
mousePrP[Y169A,Y225A,Y226A] & (R164), (R229), (R230), (K185, K194), (R136, K220), (R148), (R151)\\ \hline  
elkPrP & (R164), (R228), (K185), (HHT121), (K194), (R156), (R136, R220), (R148), (R151)\\ \hline
pigPrP & (R164), (R228), (K185, H187, K194, R156, R136), (K220), (R148, R151)\\ \hline
BankVolePrP      & (R164), (R229), (K185), (K194), (R156), (R148), (R151, R136), (H140), (K220)\\ \hline
TammarWallabyPrP & (R164), (R227), (K185, H187), (K194), (R156), (R148), (R136, H140, R148, R151, R156)\\ \hline
rabbitPrP-NMR    & (R163, R227), (K184, HHT124)-(K193, R155, R135, R150, H139, R147), (R147)\\ \hline
rabbitPrP-X-ray  & (R164, HHT126, K185)-(K194, H187, R136), (R228), (R156), (R148), (R151), (H140)\\ \hline
horsePrP         & (R164), (R228,K220,R136),(K194,H187,K185),(R156),(R148, R151),(H177),(K204),(R208),(HHT119)\\ \hline
\end{tabular} 
} \label{cliques_of_positively_charged_residues}
\end{table} 

\begin{table}[h!]
\caption{\textsf{$\pi$-$\pi$-stacking and $\pi$-cations for each of the 33 PrPs:}}
\centering
{\tiny
\begin{tabular}{|l                            |l                                              |l|} \hline
                Species                       &$\pi$-$\pi$-stacking                           &$\pi$-cations\\ \hline \hline
                mousePrP                      &F175--Y218, Y162--Y128, H187--F198 &F141--R208\\ \hline
                humanPrP                      &                                               &R164--Y169\\ \hline
                bovinePrP                     &                                               &Y128--R164\\ \hline
                SyrianHamsterPrP              &Y169--F175--Y218                               &\\ \hline
                dogPrP                        &                                               &\\ \hline
                catPrP                        &                                               &Y150--R156\\ \hline
                sheepPrP                      &F141--Y150, Y169--F175--Y218                   &\\ \hline
                mousePrP[N174T]               &F141--Y150, Y169--F175--Y218                   &Y128--R164\\ \hline
                humanPrP[Q212P]-M129          &                                               &R228--H237\\ \hline
                humanPrP-M129                 &                                               &\\ \hline
                rabbitPrP[S173N]-NMR          &                                               &L124--Y127\\ \hline
                rabbitPrP[I214V]-NMR          &H139--Y149                                     &Y148--R155\\ \hline
                rabbitPrP[S170N]-X-ray        &Y169--F175                                     &F141--R208\\ \hline
                rabbitPrP[S174N]-X-ray        &Y169--F175                                     &F141--R208\\ \hline
                rabbitPrP[S170N,S174N]-X-ray  &Y169--F175                                     &F141--R208\\ \hline \hline
                mousePrP - at 37 $\mathring{\text{}}$C             &                          &R164--Y169\\ \hline
                mousePrP[V166A]               &Y169--F175                                     &\\ \hline
                mousePrP[D167S] - at 20 $\mathring{\text{}}$C      &F175--Y218                &\\ \hline
                mousePrP[D167S,N173K]         &F175--Y218, H187--F198                         &\\ \hline
                mousePrP[Y169G]               &F141--Y150, F175--Y218, Y225--Y226             &Y128--R164\\ \hline
                mousePrP[Y169A]               &W145--Y149, H187--F198                       &\\ \hline
                mousePrP[S170N]               &Y225--Y226                                     &\\ \hline
                mousePrP[S170N,N174T]         &                                               &R164--Y169\\ \hline
                mousePrP[F175A]               &Y163--Y218                                     &F141--R208\\ \hline
                mousePrP[Y225A,Y226A]         &Y169--F175--Y218                               &\\ \hline
                mousePrP[Y169A,Y225A,Y226A] - at 20 $\mathring{\text{}}$C &                   &\\ \hline
                elkPrP                        &Y169--F175--Y218                               &\\ \hline
                pigPrP                        &                                               &\\ \hline
                BankVolePrP                   &Y169--F175--Y218                               &R164--Y169\\ \hline
                TammarWallabyPrP              &                                               &R156--F198\\ \hline
                rabbitPrP-NMR                 &F140--Y149                                     &L124--Y127\\ \hline
                rabbitPrP-X-ray               &Y169--F175                                     &\\ \hline
                horsePrP                      &                                               &R156--F198\\ \hline      
\end{tabular} 
} \label{pi_stacks}
\end{table}   

\section{A concluding remark}
In optimization, especially for solving large scale or complex or both optimization problems, the hybrid of optimization (local search or global search) methods is very necessary, and very effective and efficient for solving optimization problems. In molecular mechanics, to optimize its potential energy, even just one part of it e.g. the Lennard-Jones potential, is still a challenge to optimization methods; the hybrid idea is very helpful and useful. An application to prion protein structures is then done by the hybrid idea. Focusing on the $\beta$2-$\alpha$2 loop of prion protein structures, we found (i) the species that has the clearly and highly ordered $\beta$2-$\alpha$2 loop usually owns a 3$_{10}$-helix in this loop, (ii) a ``$\pi$-circle" Y128--F175--Y218--Y163--F175--Y169--R164--Y128(--Y162) is around the $\beta$2-$\alpha$2 loop. In conclusion, this paper proposes a hybrid idea of optimization methods to efficiently solve the potential energy minimization problem and the LJ clusters problem. We first reviewed several most commonly used classical potential energy functions and their optimization methods used for energy minimization, as well as some effective computational optimization methods used to solve the problem of Lennard-Jones clusters. In addition, we applied this hybrid idea to construct molecular structures of prion amyloid fibrils at AGAAAAGA segment, by which we provided the additional insight for the $\beta$2-$\alpha$2 loop of prion protein structures. This study should be of interest to the protein structure field. 

\section*{Acknowledgments}
This research was supported by a Victorian Life Sciences Computation Initiative (VLSCI) grant numbered VR0063 on its Peak Computing Facility at the University of Melbourne, an initiative of the Victorian Government (Australia). This paper is dedicated to Professor Alexander M. Rubinov in honour of his 75th birthday and this paper was reported in the Workshop on Continuous Optimization: Theory, Methods and Applications, 16-17 April 2015, Ballarat, Australia.


\begin{thebibliography} {99}
\bibitem{abraham_etal_GROMACS_2014}
Abraham M.J., van der Spoel D., Lindahl E., Hess B., and GROMACS development team (2014) GROMACS User Manual version 5.0.4, www.gromacs.org.

\bibitem{bagirovkm2014}
Bagirov A., Karmitsa N., Makela M.M. (2014) Introduction to Nonsmooth Optimization -Theory, Practice and Software, Springer, ISBN 978--3--319--08113--7.

\bibitem{barrongrs1999}
Barron C., Gomez S., Romero D., Saavedra A. (1999) A genetic algorithm for Lennard-Jones atomic clusters. Applied Mathematics Letters 12(7) 85--90.

\bibitem{barronromero2005}
Barron-Romero C. (2005) Minimum search space and efficient methods for structural cluster optimization. arXiv0504030v5.

\bibitem{bhandarkar_etal_NAMD_2012}
Bhandarkar M., Bhatele A., Bohm E., Brunner R., et al. (2012). NAMD User's Guide Version 2.9, University of Illinois and Beckman Institute, Urbana USA.

\bibitem{biljan_etal2012a}
Biljan, I., Giachin, G., Ilc, G., Zhukov, I., Plavec, J.,  Legname, G. (2012a). Structural basis for the protective effect of the human prion protein carrying the dominant-negative E219K polymorphism.  Biochemical Journal 446, 243--251.

\bibitem{biljan_etal2012b}
Biljan, I., Ilc, G., Giachin, G., Legname, G.,  Plavec, J. (2012b). Structural rearrangements at physiological pH: nuclear magnetic resonance insights from the V210I human prion protein mutant.  Biochemistry 51,7465--7474.

\bibitem{biljanigrzpl2011}
Biljan, I., Ilc, G., Giachin, G., Raspadori, A., Zhukov, I., Plavec, J.,  Legname, G. (2011). Toward the molecular basis of inherited prion diseases: NMR structure of the human prion protein with V210I mutation.  Journal of Molecular Biology 412, 660--673.

\bibitem{calvodw2001}
Calvo F., Doye J.P.K., Wales D.J. (2001) Quantum partition functions from classical distributions. Application to rare-gas clusters. Journal Chem. Phys. 114: 7312--29.

\bibitem{calzolai_etal2000}
Calzolai, L., Lysek, D.A., Guntert, P., Von Schroetter, C., Zahn, R., Riek, R., Wuthrich, K. (2000). NMR structures of three single-residue variants of the human prion protein.  Proceedings of the National Academy of Sciences USA 97, 8340--8345.

\bibitem{amber14}
D.A. Case, V. Babin, J.T. Berryman, R.M. Betz, Q. Cai, D.S. Cerutti, T.E. Cheatham, III, T.A. Darden, R.E. Duke, H. Gohlke, A.W. Goetz, S. Gusarov, N. Homeyer, P. Janowski, J. Kaus, I. Kolossváry, A. Kovalenko, T.S. Lee, S. LeGrand, T. Luchko, R. Luo, B. Madej, K.M. Merz, F. Paesani, D.R. Roe, A. Roitberg, C. Sagui, R. Salomon-Ferrer, G. Seabra, C.L. Simmerling, W. Smith, J. Swails, R.C. Walker, J. Wang, R.M. Wolf, X. Wu and P.A. Kollman (2014).  AMBER 14, University of California, San Francisco.

\bibitem{christen_etal2009}
Christen, B., Hornemann, S., Damberger, F.F.,  Wüthrich, K. (2009). Prion protein NMR structure from tammar wallaby (Macropus eugenii) shows that the $\beta$2-$\alpha$2 loop is modulated by long-range sequence effects.  Journal of Molecular Biology 389, 833--845.

\bibitem{christen_etal2012}
Christen, B., Hornemann, S., Damberger, F.F.,  Wüthrich, K. (2012). Prion protein mPrP [F175A](121--231): structure and stability in solution.  Journal of Molecular Biology 423, 496--502.

\bibitem{colemansw1997}
Coleman T., Shalloway D., Wu Z.J. (1997) A  parallel build-up algorithm for global energy minimizations of molecular clusters using effective energy simulated annealing. Journal Global Optimization 4: 171--85.

\bibitem{dambergercphw2011}
Damberger, F.F., Christen, B., Pérez, D.R., Hornemann, S.,  Wüthrich, K. (2011). Cellular prion protein conformation and function.  Proceedings of the National Academy of Sciences USA 108, 17308--17313.

\bibitem{deaventmh1996}
Deaven D.M.,  Tit N., Morris J.R., Ho K.M. (1996) Structural optimization of Lennard-Jones clusters by a genetic algorithm.  Chem.  Phys.  Lett. 256: 195--200.

\bibitem{doyew1995}
Doye J.P.K., Wales D.J. (1995) Magic numbers and growth sequences of small face-centred-cubic and decahedral  clusters.  Chem. Phys. Lett. 247: 339--47.

\bibitem{doyewb1995}
Doye J.P.K., Wales D.J., Berry R.S. (1995) The effect of the range of the potential on the structures of clusters. Journal Chem. Phys. 103(10): 4234--49.

\bibitem{farges_etal1985} 
Farges J.,  de Feraudy M.F., Raoult B., Torchet G. (1985) Cluster models made of double icosahedron units. Surf. Sci. 156:  370--8.

\bibitem{freemand1985}
Freeman D.L., Doll J.D. (1985) Quantum Monte Carlo study of the thermodynamic properties of argon clusters: the homogeneous nucleation of argon in argon vapor and ``magic number"  distributions in argon vapor.  Journal Chem.  Phys. 82: 462--71.

\bibitem{gossertblfw2005}
Gossert, A.D., Bonjour, S., Lysek, D.A., Fiorito, F.,  Wüthrich, K. (2005). Prion protein NMR structures of elk and of mouse/elk hybrids.  Proceedings of the National Academy of Sciences USA 102, 646--650.

\bibitem{hoarep1971a}
Hoare  M.R.,  Pal P. (1971) Statics and stability of small cluster neclei.  Nature (Physical Sciences) 230, 5--8.

\bibitem{hoarep1971b}
Hoare M.R.,  Pal P. (1971) Physical cluster mechanic: statics and energy surfaces for monatomic systems.  Advances in Chemical Physics 20, 161--96. 

\bibitem{hoarep1972}
Hoare M.R., Pal P. (1972). Geometry and stability of ``spherical"  f.c.c. microcrys- tallites.  Nature (Physical Sciences) 236, 35--7.

\bibitem{horowitzs1990}
Horowitz E., Sahni S. (1990)  Data Structures in Pascal, Computer  Science Press.

\bibitem{huangc2015}
Huang, D.,  Caflisch, A. (2015). Evolutionary conserved Tyr169 stabilizes the $\beta$2-$\alpha$2 loop of prion protein.  Journal of American Chemical Society 137, 2948--57.

\bibitem{ile_etal2010}
Ilc, G., Giachin, G., Jaremko, M., Jaremko, L., Benetti, F., Plavec, J., Zhukov, I.,  Legname, G. (2010). NMR structure of the human prion protein with the pathological Q212P mutation reveals unique structural features.  PLoS ONE 5(7), e11715.

\bibitem{kong_etal2013}
Kong, Q., Mills, J.L., Kundu, B., Li, X., Qing, L., Surewicz, K., ... Surewicz, W.K. (2013). Thermodynamic stabilization of the folded domain of prion protein inhibits prion infection in vivo.  Cell Reports 4, 248--254.

\bibitem{kurt_etal2014a}
Kurt TD, Bett C, Fernández-Borges N, Joshi-Barr S, Hornemann S, Rülicke T, ... Sigurdson CJ (2014a). Prion transmission prevented by modifying the $\beta$2-$\alpha$2 loop structure of host PrP$^\text{C}$.  Journal of Neuroscience 34, 1022--1027.

\bibitem{kurt_etal2014b}
Kurt, T.D., Jiang, L., Bett, C., Eisenberg, D.,  Sigurdson, C.J. (2014b) A proposed mechanism for the promotion of prion conversion involving a strictly conserved tyrosine residue in the $\beta$2-$\alpha$2 loop of PrP$^\text{C}$.  Journal of Biological Chemistry 289, 10660--10667.

\bibitem{laixh2011a}
Lai X.J., Xu R.C., Huang W.Q. (2011a) Prediction of the lowest energy configuration for Lennard-Jones clusters. Science China Chemistry 54(6): 985--91.

\bibitem{laixh2011b}
Lai X.J., Xu R.C., Huang W.Q. (2011b) Geometry optimization of bimetallic clusters using an efficient heuristic method. Journal Chem. Phys. 135: 164109.

\bibitem{leary1997}
Leary R.H. (1997) Global  optima of Lennard-Jones clusters. Journal of Global Optimization 11: 35--53.

\bibitem{leeahkssy2010}
Lee, S., Antony, L., Hartmann, R., Knaus, K.J., Surewicz, K., Surewicz, W.K., Yee, V.C. (2010). Conformational diversity in prion protein variants influences intermolecular beta-sheet formation.  EMBO Journal 29, 251--62.

\bibitem{locatellis2008}
Locatelli M., Schoen F. (2008) Structure prediction and global optimization. Optima Mathematical Programming Society Newsletter USA 76,1--8.

\bibitem{mackay1962}
Mackay A.L. (1962)  A  dense  non-crystallographic  packing  of equal  spheres.  Acta Crystallographica 15: 916--8.

\bibitem{northby1987}
Northby J.A. (1987) Structure and binding  of the Lennard-Jones clusters:  13 = N =147. Journal Chem.  Phys.  87  (10): 6166--77.

\bibitem{paquetv2015}
Paquet E., Viktor H.L. (2015) Molecular dynamics, Monte Carlo simulations, and Langevin dynamics: a computational review. BioMed Research International, DOI: 10.1155/2015/183918.

\bibitem{perezdw2010}
Perez, D.R., Damberger, F.F.,  Wuthrich, K. (2010). Horse prion protein NMR structure and comparisons with related variants of the mouse prion protein.  Journal of Molecular Biology 400, 121--128.

\bibitem{perezw2008}
Perez, D.R.,  Wuthrich, K. (2008). NMR structure of the bank vole prion protein at 20 degrees C contains a structured loop of residues 165-171.  Journal of Molecular Biology 383, 306--312.

\bibitem{romerobg1999}
Romero D., Barron C., Gomez S. (1999) The optimal geometry of Lennard-Jones clusters: 148--309. Computer Physics Communications 123(1-3): 87--96.

\bibitem{stote_etal1999}
Stote R., Dejaegere A., Kuznetsov D., Falquet L. (1999). Molecular dynamics simulations CHARMM. 
http://www.ch.embnet.org/MD\_tutorial/

\bibitem{sunz2001}
Sun J., Zhang J.P. (2001) Global Convergence of Conjugate Gradient Methods without Line Search. Annals of Operations Research 103(1--4): 161--73.

\bibitem{sweeting_etal2013}
Sweeting, B., Brown, E., Khan, M.Q., Chakrabartty, A.,  Pai, E.F. (2013). N-terminal helix-cap in $\alpha$-helix 2 modulates $\beta$-state misfolding in rabbit and hamster prion proteins.  PLoS One 8(5), e63047.

\bibitem{takeuchi2006}
Takeuchi H. (2006) Clever and efficient method for searching optimal geometries of Lennard- Jones clusters. Journal of Chemical Information and Modeling 46(5): 2066--70.

\bibitem{walesd1997} 
Wales D.J., Doye J.P.K. (1997) Global optimization by Basin-Hopping and the lowest energy structures of Lennard-Jones clusters containing up to 110 atoms. Journal Phys. Chem. A 101: 5111--6.

\bibitem{wenlxpyhl2010}
Wen, Y., Li, J., Xiong, M.Q., Peng, Y., Yao, W.M., Hong, J.,  Lin, D.H. (2010). Solution structure and dynamics of the I214V mutant of the rabbit prion protein.  PLoS One 5(10), e13273.

\bibitem{wenlyxpxl2010}
Wen, Y., Li, J., Yao, W.M., Xiong, M.Q., Hong, J., Peng, Y., ... Lin, D.H. (2010). Unique structural characteristics of the rabbit prion protein.  Journal of Biological Chemistry 285, 31682--31693.

\bibitem{wille1987}
Wille L.T. (1987) Minimum-energy configurations of atomic clusters: new results obtained by simulated annealing. Chem. Phys. Lett. 133: 405-.

\bibitem{wolf1998}
Wolf M. (1998) Journal Phys. Chem. A. 102(30).

\bibitem{xiangccs2004}
Xiang Y., Cheng L., Cai W., Shao X. (2004b) Structural distribution of Lennard-Jones clusters containing 562 to 1000 atoms. Journal Phys. Chem. A 108: 9516.

\bibitem{xiangjcs2004}
Xiang Y.,Jiang H., Cai W., Shao X. (2004a) An efficient method based on lattice construction and the genetic algorithm for optimization of large Lennard-Jones clusters. Journal Phys. Chem. A 108: 3586--92.

\bibitem{xue1994a}
Xue G.L. (1994a) Molecular conformation on the CM-5 by parallel two-level simulated annealing. Journal of Global Optimization 4: 187--208.

\bibitem{xue1994b}
Xue G.L. (1994b) Improvement on the Northby algorithm for molecular  confirmation: better solutions. Journal of Global Optimization 4: 425--40.

\bibitem{yexh2011}
Ye T.,  Xu R.C., Huang W.Q. (2011) Global optimization of binary Lennard-Jones clusters using three perturbation operators. Journal Chem. Inf. Model. 51 (3): 572–-7.

\bibitem{zahngvsw2003}
Zahn, R., Guntert, P., von Schroetter, C.,  Wüthrich, K. (2003). NMR structure of a variant human prion protein with two disulfide bridges.  Journal of Molecular Biology 326, 225--234.

\bibitem{zhang2010}
Zhang J.P. (2010) Studies on the structural stability of rabbit prion probed by molecular dynamics simulations of its wild-type and mutants. Journal of Theoretical Biology 264(1): 119--22.

\bibitem{zhang2011a}
Zhang J.P. (2011a) Optimal molecular structures of prion AGAAAAGA amyloid fibrils formatted by simulated annealing. Journal of Molecular Modeling 17(1): 173--9. 

\bibitem{zhang2011b}
Zhang J.P. (2011b) Practical Global Optimization Computing Methods in Molecular Modelling - for Atomic-resolution Structures of Amyloid Fibrils, LAP LAMBERT Academic Publishing, ISBN 978--3--8465--2139--7.

\bibitem{zhang2011c}
Zhang J.P. (2011c) An effective simulated annealing refined replica exchange Markov chain Monte Carlo algorithm for the infectious disease model of H1N1 influenza pandemic. World Journal of Modelling and Simulation 7(1) 29–39.

\bibitem{zhang2011d}
Zhang J.P. (2011d) The structural stability of wild-type horse prion protein. Journal of Biomolecular Structure and Dynamics 29(2) : 369--77.

\bibitem{zhang2011e}
Zhang J.P. (2011e) Comparison studies of the structural stability of rabbit prion protein with human and mouse prion proteins. Journal of Theoretical Biology 269(1): 88--95.

\bibitem{zhang2011f}
Zhang J.P. (2011f) Molecular dynamics – practical application – mechanism underlying the resistance to prion diseases in rabbits: a technology review from MIT and in 7 languages, LAP LAMBERT Academic Publishing, ISBN 978-3-8465-4843-1.

\bibitem{zhang2014}
Zhang J.P. (2014) Simulated annealing: in mathematical global optimization computation, hybrid with local or global search, and practical applications in crystallography and molecular modelling of prion amyloid fibrils.[Simulated Annealing: Strategies, Potential Uses and Advantages, Editors Prof. Dr. Marcos Tsuzuki \& Prof. Dr. Thiago de Castro Martins, NOVA Science Publishers, ISBN 978-1-63117-268-7] Chapter 1:1-34. 

\bibitem{zhang2015b}
Zhang J.P. (2015) A survey on $\pi$-$\pi$ stackings and $\pi$-cations in prion protein structures $\sim$ A `Quick Reference Card'. Biochem Pharmacol (Los Angel) 4(3):e175, doi: 10.4172/2167-0501.1000e175. 

\bibitem{zhang2015a}
Zhang J.P., Wang F., Chatterjee S. (2015) Molecular dynamics studies on buffalo prion protein. Journal of Biomolecular Structure and Dynamics, doi: 10.1080/07391102.2015.1052849. 

\bibitem{zhanghwwz2012}
Zhang, J.P., Hou Y.T., Wang Y.J., Wang C.Y., Zhang X.S. (2012) The LBFGS quasi-Newtonian method for molecular modeling prion AGAAAAGA amyloid fibrils. Natural Science 4(12A): 1097--108.

\bibitem{zhanggy2011}
Zhang J.P., Gao D.Y., Yearwood J. (2011) A novel canonical dual computational approach for prion AGAAAAGA amyloid fibril molecular modeling. Journal of Theoretical Biology 284(1): 149--57.

\bibitem{zhangl2011}
Zhang J.P., Liu D.D.W. (2011) Molecular dynamics studies on the structural stability of wild-type dog prion protein. J of Biomolecular Structure and Dynamics 28(6): 861--9.

\bibitem{zhangsw2011}
Zhang J.P., Sun J., Wu C.Z. (2011) Optimal atomic-resolution structures of prion AGAAAAGA amyloid fibrils. Journal of Theoretical Biology 279(1): 17--28.

\bibitem{zhangszss2000}
Zhang, Y., Swietnicki, W., Zagorski, M.G., Surewicz, W.K.,  Soennichsen, F.D. (2000). Solution structure of the E200K variant of human prion protein. Implications for the mechanism of pathogenesis in familial prion diseases.   Journal of Biological Chemistry 275, 33650--33654.

\bibitem{zhangwz2015}
Zhang J.P., Wang F. Zhang Y.L. (2015) Molecular dynamics studies on the NMR structures of rabbit prion protein wild-type and mutants: surface electrostatic charge distributions. J Biomol Struct Dyn 33(6): 1326-35.

\bibitem{zhangz2013}
Zhang, J.P., Zhang, Y.L. (2013). Molecular dynamics studies on 3D structures of the hydrophobic region PrP(109--136).  Acta Biochimicaet Biophysica Sinica 45, 509--519.

\bibitem{zhangz2014}
Zhang, J.P., Zhang, Y.L. (2014). Molecular dynamics studies on the NMR and X-ray structures of rabbit prion proteins.  Journal of Theoretical Biology 342, 70--82.

\bibitem{zhao_etal2012}
Zhao, H., Liu, L.L., Du, S.H., Wang, S.Q.,  Zhang, Y.P. (2012). Comparative analysis of the Shadoo gene between cattle and buffalo reveals significant differences.  PLoS One 7(10), e46601.

\bibitem{zhong2010}
Zhong, L.H. (2010). Exposure of hydrophobic core in human prion protein pathogenic mutant H187R.  Journal of Biomolecular Structure and Dynamics 28, 355--361. 
\end{thebibliography}
\end{document}